\documentclass{opt2022} 

\usepackage{capt-of}
\usepackage{float}
\usepackage{wrapfig}
\usepackage{algorithmic}
\usepackage{algorithm}
\usepackage{multirow}

\newcommand{\ayan}[1]{{\bf \color{purple} #1}}

\newcommand{\R}{\mathbb{R}}
\newcommand{\T}{\top}

\newcommand{\lam}{\lambda}
\newcommand{\dt}{\delta}
\newcommand{\Dt}{\Delta}

\newcommand{\al}{\alpha}

\newcommand{\ep}{\epsilon}

\newcommand{\argmax}{\mathop{{\rm argmax}}}

\newcommand{\td}{\tilde}

\newcommand{\Ex}{\mathbb{E}}

\newcommand{\cI}{\mathcal{I}}

\newcommand{\cS}{\mathcal{S}}
\newcommand{\cT}{\mathcal{T}}

\newcommand{\cX}{\mathcal{X}}

\newcommand{\ext}{{\rm ext}}
\newcommand{\var}{{\rm var}}

\title[Fast LP Solver for Personalized Recommendation with Diversity Constraints]{A Light-speed Linear Program Solver for Personalized Recommendation with Diversity Constraints}

\optauthor{%
\Name{Haoyue Wang}\thanks{This work was done when Haoyue Wang was an intern at LinkedIn during summer 2022.} \Email{wanghaoyue3@outlook.com}\\
\addr LinkedIn, Sunnyvale, CA \\
(Massachusetts Institute of Technology, Cambridge, MA.)
\AND
\Name{Miao Cheng} \Email{miacheng@linkedin.com}\\
\addr LinkedIn, Sunnyvale, CA
\AND
\Name{Kinjal Basu} \Email{kbasu@linkedin.com}\\
\addr LinkedIn, Sunnyvale, CA
\AND
\Name{Aman Gupta} \Email{amagupta@linkedin.com}\\
\addr LinkedIn, Sunnyvale, CA
\AND
\Name{Keerthi Selvaraj} \Email{keselvaraj@linkedin.com}\\
\addr LinkedIn, Sunnyvale, CA
\AND
\Name{Rahul Mazumder}\thanks{This work was done when 
Rahul Mazumder was a consultant for LinkedIn.} \Email{rmazumder@linkedin.com}\\
\addr LinkedIn, Sunnyvale, CA \\
(Massachusetts Institute of Technology, Cambridge, MA.)
}

\begin{document}

\maketitle

\begin{abstract}%
We study a structured linear program (LP) that emerges in the need of ranking candidates or items in personalized recommender systems. 
Since the candidate set is only known in real time, the LP also needs to be formed and solved in real time. Latency and user experience are major considerations, requiring the LP to be solved within just a few milliseconds. Although typical instances of the problem are not very large in size, this stringent time limit appears to be beyond the capability of most existing (commercial) LP solvers, which can take $20$ milliseconds or more to find a solution. 
Thus, reliable methods that address the real-world complication of latency become necessary.
In this paper, we propose a fast specialized LP solver for a structured problem with diversity constraints. 
Our method solves the dual problem, making use of the piece-wise affine structure of the dual objective function, with an additional screening technique that helps reduce the dimensionality of the problem as the algorithm progresses. Experiments reveal that our method can solve the problem within roughly 1 millisecond, yielding a 20x improvement in speed over efficient off-the-shelf LP solvers. This speed-up can help improve the quality of recommendations without affecting user experience, highlighting how optimization can provide solid orthogonal value to machine-learned recommender systems. 
\end{abstract}


    \section{Introduction}
	
Linear programs (LPs) arise naturally in several applications across various scientific and industrial disciplines. Examples include scheduling \cite{hanssmann1960linear}, matching entities \cite{azevedo2016matching,zheng2017online,borisyuk2017lijar}, network flow \cite{ahuja1988network,bazaraa2008linear}, recommender systems \cite{linden2003amazon,agarwal2012personalized,agarwal2014activity,agarwal2015personalizing} and many others \cite{gupta2016email,gupta2017optimizing}. Despite the significant development of LP solvers after decades of research, modern applications bring new challenges in solving LP problems, either with extreme-scale datasets \cite{basu2020eclipse}, or with strict restrictions on time budget (e.g. in an online setting where LPs need to be solved in real-time). These challenges have sparked recent explorations for more efficient algorithms to solve LPs in both academia and industry, see e.g. \cite{basu2020eclipse,lin2021admm,applegate2021faster,applegate2021practical,song2021coordinate, DBLP:journals/corr/abs-2103-05277}.

In this paper, we focus on a framework for re-ranking  candidates in recommender systems under diversity or grouping constraints \cite{carbonell1998use}. 
Consider a social media platform that makes personalized recommendations for its users. 
Suppose there are $m$ candidates (which could be items \cite{agarwal2015personalizing} or people \cite{zheng2017online}) for generating recommendations for user $u$; each candidate $i$ has a utility score $c_i \in \R$ that is generated by a machine learning model, representing the probability that $u$ may interact with this candidate. Suppose there are $n$ slots to present the recommendations, with weight $w_j>0$ representing the importance of the slot $j$. In particular, we assume $w_1>w_2 > \cdots > w_n$, i.e., slots in the `front' are more important. In light of the applications we consider, we assume that $n$ is much smaller than $m$: for example, $n$ equals $10\sim 20$, while $m$ can range from hundreds to thousands. 
The recommendations and scores are specific to a given user $u$, and vary across users. 
	
The primary goal of the system is to find an assignment of the $m$ candidates into the $n$ slots so as to maximize the total weighted utility score. We formulate the LP problem per user $u$, and solve the LPs for all different users in parallel. Let matrix $X^u \in \R^{m\times n}$ represent a randomized\footnote{We use the term randomized to refer to the fact that the entries of $X$ can be fractional} assignment for each user $u$'s recommendations, where each $X^u_{ij}\in [0,1]$ stands for the probability that candidate $i$ is assigned to slot $j$. For simplicity, we remove the superscript of $u$ in all the following notations. Then for each user $u$, the recommendation problem can be formulated as:
	\begin{equation}\label{LP-without-diversity}
			\max_{X\in \R^{m\times n}} ~ c^\T X w
			\quad ~
			{\rm s.t.} ~  X \in \cS_{m,n}
	\end{equation}
where $c = (c_1,c_2,...,c_m)^\T \in \R^m$, $w = (w_1,w_2,...,w_n)^\T \in \R^n$, and
the set $\cS_{m,n}$ is defined by 
\begin{equation}
	\cS_{m,n} := \{ X\in \R^{m\times n}~|~ X 1_{n} \le 1_m, ~ 1_m^\T X = 1_n^\T , ~ X\ge 0 \}
\end{equation}
where $1_k$ is the vector in $\R^k$ with all coordinates being $1$ (for any $k\ge 1$). The constraint $X 1_{n} \le 1_m $ implies that for each candidate $i$, the probability that $i$ is assigned to some slot is at most $1$. The constraint $1_m^\T X = 1_n^\T $ implies that each slot must be assigned a candidate. 
Note that (\ref{LP-without-diversity}) has a simple closed-form solution: let $\sigma^*$ be a permutation on $[m]$ such that $c_{\sigma^*(1)} \ge c_{\sigma^*(2)} \ge \cdots \ge c_{\sigma^*(m)} $, then matrix $X^*$ with $X^*_{\sigma^*(j), j} = 1$ for all $j\in [n]$ and all other entries being $0$, is an optimal solution of (\ref{LP-without-diversity}). In other words, the candidates with the top $n$ scores are placed in the top $n$ slots, in descending order of the scores. 

In the system described above, only utility scores are used to make assignments. Recommendations obtained this way often lead to selection bias and can favor a particular group in terms of the overall ranking. For example, different kinds of people, different activity levels across the platform, different types of content etc. may get a very different exposure or placement in the rankings. Average ranks of one group may be much higher than those of another group when the system is solely optimized for overall utility which increases homogeneity over time.

To allow for diverse recommendations, we introduce an additional \textit{diversity constraint} to this problem. Let $a_i\in \R$ denote the diversity feature of the candidate $i$. The following LP places a constraint on the weighted sum of the diversity features:
\begin{equation}\label{LP-with-diversity}
	\max_{X\in \R^{m\times n}} ~ c^\T X w
	\quad ~
	{\rm s.t.} ~  X \in \cS_{m,n}, ~ 
	b_1 \le a^\T X w \le b_2  . 
\end{equation}
where $a = (a_1,a_2,...,a_m)^\T \in \R^m$, and $b_1$ and $b_2$ are chosen lower and upper bounds on the total weighted diversity. This allows us to incorporate a family of desired diversity restrictions with appropriate choices on $a, b_1$ and $b_2$. 
{For example, in ``People You May Know" (PYMK) recommendations, if the user (viewer) is a software engineer, then it is likely that the system would give candidates who are also software engineers a favorably higher score $c_i$. It may be better to blend more candidates who are not software engineers but share other common similarities with the user. In this case, we model $a_i=1$ if the candidate has similar job as the user, and $a_i=-1$ if not. By properly setting values $b_1$ and $b_2$ we can diversify the PYMK recommendations.}
A few more concrete examples of real-world use cases which can benefit from our specialized solver are discussed in the appendix (Section~\ref{sect:real-world-use-case}).

Unlike \eqref{LP-without-diversity}, the LP in \eqref{LP-with-diversity} does not have a closed-form solution. As an example: A leading commercial solver Gurobi can solve Problem~\eqref{LP-with-diversity} with small or medium size (e.g., $m=100$, $n=10$)
in roughly 20 milli-seconds using their proprietary solver. 
However, for an online application where data $c_i$ (and the set of candidates themselves) are frequently updated, the personalized LP problem \eqref{LP-with-diversity} needs to be solved in an ``on-click" manner: when a user clicks a button to load the recommendations, the LP problem under the hood needs to be solved in an instant, i.e. ideally within $3\sim 5$ milli-seconds. Our goal here is to develop a novel algorithm to be able to solve this LP under such strict time-requirements. 

In this paper, we introduce a specialized algorithm for problem \eqref{LP-with-diversity} that can solve an instance with $m=100$ and $n=10$ within one milli-second, which is a 20x speedup compared to commercial LP solvers. 
Our method introduces a dual variable for the diversity constraint in \eqref{LP-with-diversity} and solves the one-dimensional dual problem using a bisection method. When we are sufficiently close to an optimal dual solution, we trace the piece-wise linear structure of the dual objective function to locate the exact value of the dual optimal solution. The primal optimal solution can then be recovered from the dual optimal solution. We also introduce a screening technique that can detect and drop candidates that would not be assigned to any slots in the optimal solution. This screening technique can significantly accelerate the algorithm when the number of candidates $m$ is large. 

The remainder of the paper is organized as follows. In Section~\ref{section: main alg}, we simplify the program \eqref{LP-with-diversity} and introduce our main algorithmic framework. In Section~\ref{section: screening} we introduce the screening technique to further accelerate the algorithm. In Section~\ref{section: exp}, we present numerical results of our algorithm. We finally wrap up with a discussion in Section~\ref{section: conclusion}. Proofs of all results are placed in the appendix for easier reading of the main paper.

\section{Dual problem and our solution approach}\label{section: main alg}


\noindent 
{\bf Reduction to a one-sided diversity constraint:} In problem \eqref{LP-with-diversity}, the diversity constraint places both a lower bound $b_1$ and an upper bound $b_2$ on $a^\T X w$. However, to find the optimal solution of \eqref{LP-with-diversity}, either the lower bound or the upper bound can be dropped, as shown by the following lemma. Define $\cX_0$ to be the set of optimal solution without the diversity constraint, i.e., $	\cX_0 := \argmax \{  c^\T X w ~|~ X\in \cS_{m,n} \} $. 
\begin{lemma}\label{lemma: remove-one-constraint}
	(i) If $\min_{X\in \cX_0} \{a^\T X w\} > b_2$, then the set of optimal solutions of \eqref{LP-with-diversity} is the same as 
	\begin{equation}\label{LP-with-diversity-ub}
		\max\nolimits_{X\in \R^{m\times n}} ~ c^\T X w
		\quad ~
		{\rm s.t.} ~~~~  X \in \cS_{m,n}, ~ 
		 a^\T X w \le b_2  . 
	\end{equation}

(ii) If $\max_{X\in \cX_0} \{a^\T X w\} < b_1$, then the set of optimal solutions of \eqref{LP-with-diversity} is the same as 
\begin{equation}\label{LP-with-diversity-lb}
	\max\nolimits_{X\in \R^{m\times n}} ~ c^\T X w
	\quad ~
	{\rm s.t.} ~~~~  X \in \cS_{m,n}, ~ 
	a^\T X w \ge b_1  . 
\end{equation}

(iii) If there exists $\bar X \in \cX_0$ such that $b_1 \le a^\T \bar X w \le b_2$, then $\bar X$ is an optimal solution of \eqref{LP-with-diversity}. 

\end{lemma}

Lemma~\ref{lemma: remove-one-constraint} provides a procedure to drop one side of the diversity constraint without changing the optimal solution. To use it in practice, we need to evaluate $\cX_0$, which can be efficiently computed 
via a ranking of the coordinates of $c$. When there is no tie in the ranking of the coordinates of $c$, the optimal set $\cX_0$ will only contain one point, say $\cX_0 = \{\bar X\}$.
%

\medskip 

\noindent {\bf Properties of the dual objective function:}
Having made use of the reduction result of Lemma~\ref{lemma: remove-one-constraint}, we now discuss how to solve problem \eqref{LP-with-diversity-ub}. 
Note that problem \eqref{LP-with-diversity-lb} can also be written into a similar form as \eqref{LP-with-diversity-ub}, with $a$ replaced by $-a$ and $b_2 $ replaced by $-b_1$. To solve problem \eqref{LP-with-diversity-ub}, we introduce a dual variable $\lam \ge 0$ for the constraint $ a^\T X^* w \le b_2$, and define the dual objective function:
\begin{equation}\label{def:g}
	g(\lam) := \max\nolimits_{X\in \cS_{m,n}} c^\T Xw - \lam ( a^\T Xw - b_2 )  =
	 \max\nolimits_{X\in \cS_{m,n}}  (c-\lam a)^\T Xw + b_2 \lam 
\end{equation}
then $g$ is a piece-wise linear convex function, and the dual problem of \eqref{LP-with-diversity-ub} is
\begin{equation}\label{dual-problem}
	\min\nolimits_{\lam \ge 0} ~~ g(\lam).
\end{equation}
For each fixed $\lam \ge 0$, the value $g(\lam)$ can be evaluated by solving the maximization in \eqref{def:g}, which can be computed in closed form via
a sorting of the coordinates of $c - \lam a$. Let $\cX_{\lam}$ be the set of optimal solutions of \eqref{def:g}. Then the extreme points of $\cX_\lam$ are given by:
\begin{equation*}
	\begin{aligned}
			\ext(\cX_{\lam}) = \{ X^{(\sigma)} ~|~ \sigma &\text{ is a permutation on } [m] \text{ with } \\
			&(c-\lam a)_{\sigma(1)} \ge (c-\lam a)_{\sigma(2)} \ge \cdots \ge (c-\lam a)_{\sigma(n)} \}
	\end{aligned}
\end{equation*}
where $X^{(\sigma)} \in \cS_{m,n}$ is the matrix with $X_{\sigma(i), i } = 1$ for $i \in [n]$ and all other entries being $0$. 
If there are no ties among the top $n$ values of $c-\lam a$, the set $\ext(\cX_{\lam})$ (and also $\cX_\lam$)  contains only one element, i.e., \eqref{def:g} has a unique solution. 
Once $\cX_{\lam}$ is computed, the derivatives of $g$ can also be evaluated, via a standard application of Danskin's theorem \cite{bonnans1998optimization}. 
\begin{lemma}
	 For any $\lam>0$, the left derivative and right derivative of $g$ are given by
	\begin{equation}\label{grad-g}
		g_-'(\lam) = b_2 - \max_{X\in \cX_\lam} \{ a^\T X w \}   , \quad \text{and} \quad 
		g_+'(\lam) = b_2 - \min_{X\in \cX_\lam} \{ a^\T X w \} .
	\end{equation} 
\end{lemma}
Note that $g$ is a piece-wise linear function. For $\lam$ in the interior of a linear block, the left and right derivatives are the same: $g_-'(\lam) = g_+'(\lam)$. 
For $\lam$ at the intersection of two linear blocks (called a \textit{kink} of $g$), it holds $g_-'(\lam) < g_+'(\lam)$. 
In particular, if at some $\lam>0$, the top $n$ coordinates of $c-\lam a$ is unique, then $\cX_{\lam}$ contains only one point, and by \eqref{grad-g}, we know that $\lam$ must lie in the interior of a linear block. Furthermore, the top $n$ values of $c-\lam a$ are the same for all other $\lam$ in the interior of the same linear block. This motivates a method to trace the piece-wise linear structure of $g$ via observing the top $n$ coordinates of $c-\lam a$ as $\lam$ changes. 

For $\lam>0$, 
let $K_R(\lam)$ be the smallest kink of $g$ that is larger than $\lam$, and $K_L(\lam)$ be the largest kink of $g$ that is smaller than $\lam$. For any vector $y\in \R^m$, let $\cT_n(y)$ be the set of top $n$ coordinates of $y$ (including ties), that is, 
	\begin{equation} 
		\cT_n(y) := \Big\{ i\in [m] ~\Big|~ \# \{j\in [m] | y_j < y_i \} \le n-1 \Big\}. 
	\end{equation}
For a given $\lam>0$, let $z:= c- \lam a$, then $K_R(\lam)$ and $K_L(\lam)$ can be computed in closed form:
\begin{equation}\label{compute-KR}
	K_R(\lam) = \lam+ \min_{(i,j)} \Big\{
	\frac{z_{i} - z_j}{a_i - a_j} ~\Big|~  j\in \cT_n(z), \ i\in [m], \ a_i - a_j >0, \ z_i - z_j>0
	\Big\}
\end{equation}
\begin{equation}\label{compute-KL}
	K_L(\lam) = \lam- \min_{(i,j)} \Big\{
	\frac{z_{i} - z_j}{a_j - a_i} ~\Big|~  j\in \cT_n(z), \ i\in [m] , \ a_j - a_i >0, \ 
	z_i - z_j >0
	\Big\}.
\end{equation}
As shown in \eqref{compute-KL} and \eqref{compute-KR}, one can trace the exact value of the kinks of $g$ starting from an arbitrary $\lam>0$. As long as $|\cT_n(y)| = O(n)$, 
the computation cost of $K_R(\lam) $ and $K_L(\lam) $ is $O(mn)$ operations. 

\medskip 

\noindent {\bf Bisection and tracing for the dual problem:} To solve problem \eqref{LP-with-diversity-ub}, we solve the one-dimensional dual problem \eqref{dual-problem} using a bisection method, which maintains and updates an interval $[\lam_{\min}, \lam_{\max}]$ that contains the dual optimal solution $\lam^*$ (for simplicity, assume it is unique). In each iteration, the algorithm computes the value and derivatives of $g$ at the midpoint $\lam= (\lam_{\min} + \lam_{\max})/2$ of this interval, and consider the following cases: (1) If $g'_+(\lam)<0$, then we know $\lam  < \lam^*$, and update the interval to $[\lam, \lam_{\max}]$; (2) if $g'_-(\lam)>0$, 
then we know $\lam  > \lam^*$, and update the interval to $[\lam_{\min}, \lam]$; (3) otherwise it holds $g'_-(\lam) \le 0 \le g'_+(\lam)$, which indicates that $\lam = \lam^*$.

In the above, the bisection method shrinks the interval by a half in each iteration, so the value of {the} dual optimal solution can be quickly located in a small interval after a few iterations. When the interval is small enough, the algorithm computes a nearby kink following \eqref{compute-KL} or \eqref{compute-KR}, and checks if this kink is the exact optimal dual solution (using the signs of left and right derivatives). If the kink is the optimal solution, the algorithm is terminated immediately; otherwise we continue the bisection iterations. Details of our algorithm are summarized in Algorithm~\ref{alg: bisection} in the appendix.
Once the dual problem is solved, the primal solution can be recovered in closed form -- see Section~\ref{section: recovering primal solutions} in the appendix for details.

Note that in this algorithm, the major algebraic operations are conducted in the computations of $\cX_{\lam}$ and $K_L(\lam)$ (or $K_R(\lam)$). 
In each iteration of Algorithm~\ref{alg: bisection}, to compute $\cX_{\lam}$, one needs to sort the vector $c - \lam a$, which takes $m\log(m)$ operations. 
When $\lam_{\max} - \lam_{\min} < \Dt$, an additional tracing step is needed to compute $K_L(\lam) $ or $K_R(\lam)$, which takes $O(mn)$ operations. Therefore, the cost of Algorithm~\ref{alg: bisection} scales (at least) linearly with respect to $m$ which can get expensive when $m$ is large. In Section~\ref{section: screening} we discuss a screening procedure to accelerate the algorithm in such a case.


\section{A screening procedure when $m $ is large}\label{section: screening}

When the number of candidates $m$ is large, 
we propose a screening technique that is used in the process of Algorithm~\ref{alg: bisection}, which can quickly detect and discard candidates that will not be in the optimal solution.
Note that the candidates that are present in the optimal solution must be the top $n$ coordinates of the vector $c - \lam^* a$, where $\lam^*$ is the dual optimal solution. Although $\lam^*$ is not available at the beginning, after a few iterations of Algorithm~\ref{alg: bisection}, we have an interval $[ \lam_{\min} , \lam_{\max}]$ that can approximately locate the value of $\lam^*$. If at some point $\lam$ that is close to $\lam^*$, we observe that some candidate $i$ has a very low rank (i.e., far from being in the top $n$) in the sorted sequence $c-  \lam a$, it is a strong evidence that $i$ is not in the top $n$ elements of vector $c-  \lam^* a$. The following lemma provides a precise statement underlying this intuition. 

\begin{lemma}\label{lemma: screening}
	Let $\lam^*$ be the solution of \eqref{dual-problem}, and let $\lam_{\min}$, $\lam_{\max}$ satisfy $ \lam_{\min} \le \lam^* \le \lam_{\max} $. For any $i\in [m]$, if there exist indices 
	$i_1,....,i_n\in [m]$ such that
	\begin{equation}\label{ineq-min-max}
		(c-\lam_{\min} a)_i < \min_{ k\in [n] } \{ (c-\lam_{\min} a)_{i_k} \}, \quad \text{and} \quad 
		(c-\lam_{\max} a)_i < \min_{ k\in [n] } \{ (c-\lam_{\max} a)_{i_k} \}
	\end{equation}
	Then the $i$-th row of the solution of \eqref{LP-with-diversity-ub} is zero. 
\end{lemma}

According to Lemma~\ref{lemma: screening}, to prove that candidate $i \in [m]$ is not selected in the optimal assignment, it suffices to find a set of ``good" indices $i_1,....,i_n$ such that \eqref{ineq-min-max} holds true. 
Different choices of $i_1,....,i_n$ may lead to different screening qualities. When the interval $[\lam_{\min}, \lam_{\max}]$ is small enough, if we take $i_1,....,i_n$ to be the top $n$ coordinates of $c - \lam a$ for some $\lam \in [\lam_{\min}, \lam_{\max}]$, then they may also have a high rank in the sortings of both $ c - \lam_{\min} a$ and $c - \lam_{\max} a$, and can be used as good candidates for the screening in \eqref{ineq-min-max}. 

We apply a screening step in each iteration of Algorithm~\ref{alg: bisection}. In particular, we take $i_1, ...., i_n$ as the top $n$ coordinates of $c - \bar \lam a$ for $\lam = (\lam_{\min} + \lam_{\max})/2$, and check the inequalities in \eqref{ineq-min-max} for all other $i \in [m]$. If for some $i$ the inequalities \eqref{ineq-min-max} hold true, then $i$ can be discarded, and the dimension of the problem is reduced. This is summarized in Algorithm~\ref{alg: bisection+screening} in appendix.

    \section{Numerical experiments}\label{section: exp}

We conduct experiments on synthetic and real datasets to {assess the performance of our algorithm}.
{All the experiments are conducted on a Macbook Pro with 4 CPUs and 16GB RAM}. 

\medskip

\noindent {\bf Speedup by the screening procedure}. 
To show the advantage of the screening steps introduced in Section~\ref{section: screening}, we compare the runtimes of Algorithms~\ref{alg: bisection} and~\ref{alg: bisection+screening} on synthetic examples with $n = 10$ and different values of $m$. See Section~\ref{section: data generation} for data generation procedure.

\begin{table}[ht]
\begin{minipage}[b]{0.45\linewidth}
\footnotesize
\begin{tabular}{|c|cc|cc|}
		\hline
		\multirow{2}{*}{} & \multicolumn{2}{c|}{$n=10$} & \multicolumn{2}{c|}{$n=30$} \\ \cline{2-5}
		& Alg 2       & Gurobi      & Alg 2       & Gurobi      \\ \hline
		$m=100$             &      0.8       &     14.8        &      1.7       &      25.4       \\ 
		$m=300$             &      1.2       &     20.3        &      2.1       &      60.1       \\ 
		$m=1000$            &      2.5       &     52.9        &      3.7       &      166.6       \\ 
		$m=3000$           &      5.6       &     115.8        &      7.3       &      526.9       \\ 
		$m=10000$           &      14.3       &     336.5        &      19.1       &      1765.3       
		\\ \hline 
	\end{tabular}
	\vspace{2mm}
	\caption{Runtime comparison of Algorithm~\ref{alg: bisection+screening} with Gurobi. Reported runtimes are in milli-seconds. }

\label{table: compare Gurobi}
\end{minipage}\hfill
\begin{minipage}[b]{0.46\linewidth}
\centering
\includegraphics[height=3.3cm,width=5.5cm]{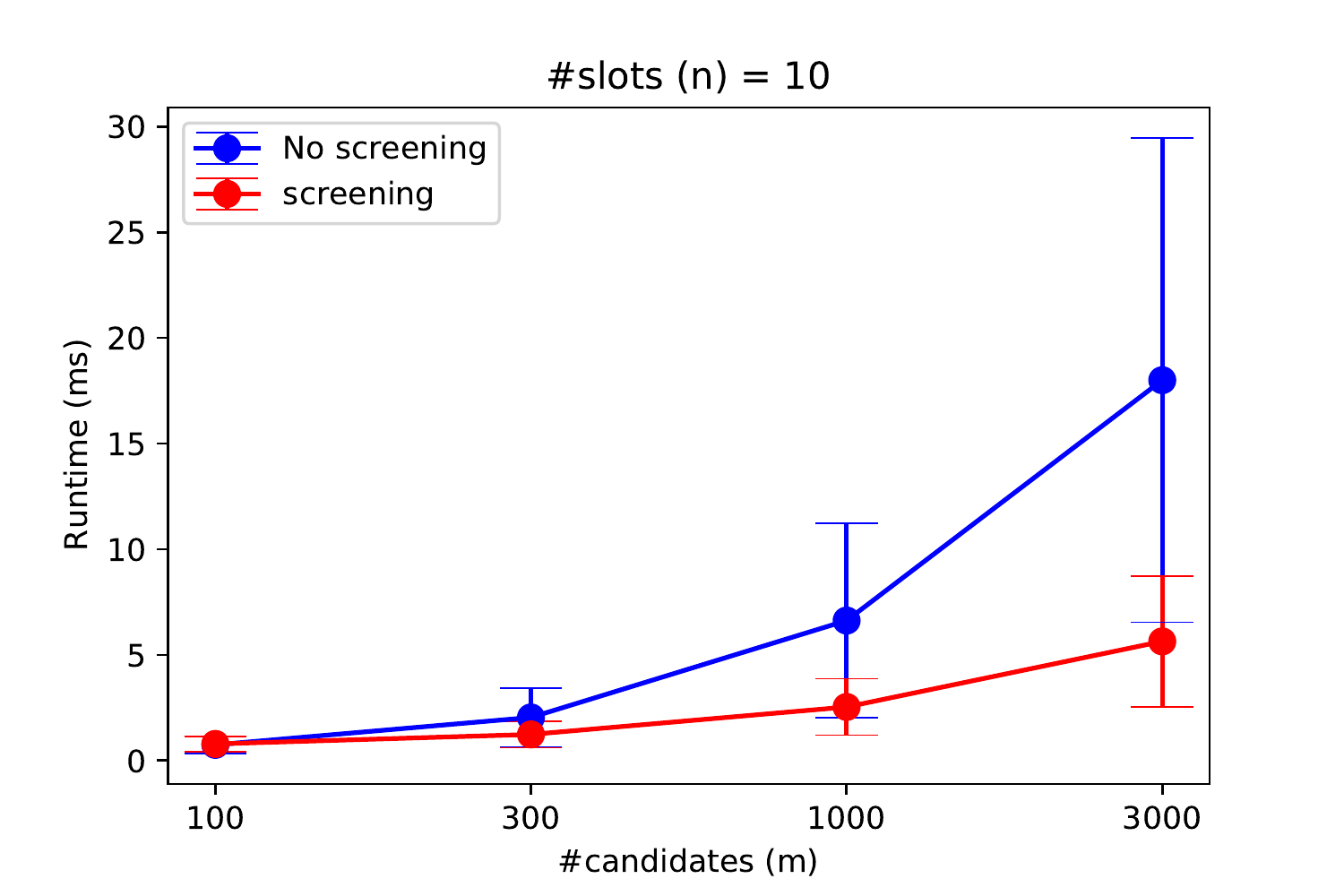}
\captionof{figure}{Runtime comparison of Algorithm~\ref{alg: bisection} (no screening) and Algorithm~\ref{alg: bisection+screening} (screening).}
	\label{figure: screening}
\end{minipage}
\end{table}


Figure~\ref{figure: screening} shows runtime comparison of Algorithms~\ref{alg: bisection} and~\ref{alg: bisection+screening} with the number of slots $n = 10$ and number of candidates $m\in \{100, 300, 1000, 3000\}$. The results are the average of 20 independent replications (with standard errors). As shown in Figure~\ref{figure: screening}, when $m=100$, the runtimes of Algorithm~\ref{alg: bisection} and Algorithm~\ref{alg: bisection+screening} are close (both are roughly one ms). For larger values of $m$, Algorithm~\ref{alg: bisection+screening} is faster than Algorithm~\ref{alg: bisection}, and the gap widens as $m$ increases. In particular, when $m=3000$, Algorithm~\ref{alg: bisection+screening} achieves a $3\sim 4$ times acceleration over Algorithm~\ref{alg: bisection}. 

\noindent
{\bf Comparison with Gurobi}. 
We compare the performance of Algorithm~\ref{alg: bisection+screening} with the commercial LP solver Gurobi~\cite{gurobi} on the synthetic data as described above. Gurobi is used under the default settings. {Note that both methods find the exact optimal solution of the problem}. 
Table~\ref{table: compare Gurobi} compares the runtimes of Algorithm~\ref{alg: bisection+screening} and Gurobi on instances of problem \eqref{LP-with-diversity} with $m\in \{100, 300, 1000, 3000, 10000\}$ and $n \in \{10,30\}$. As shown by Table~\ref{table: compare Gurobi}, when $n = 10$, Algorithm~\ref{alg: bisection+screening} is roughly $20$ times faster than Gurobi, and the acceleration is more significant for larger values of $m$. For $n = 30$, the gap between Algorithm~\ref{alg: bisection+screening} and Gurobi is even larger; in particular, for the largest example with $m = 10000$ and $n=30$, Algorithm~\ref{alg: bisection+screening} achieves more than $90$ times acceleration over Gurobi.

    \section{Conclusion}\label{section: conclusion}

We develop a highly efficient specialized linear program solver for recommender system with diversity constraints, where recommendations are made in real time with time budget in milli-seconds. Our solver is efficient in both runtime and memory requirement, with an acceleration of more than $20$ times over Gurobi. Although our primary application appears in personalized recommendation, this algorithm can be applied in other ranking applications where the problem can be formulated in the form \eqref{LP-with-diversity}.

\section{Acknowledgements}
Haoyue Wang contributed to this work while he was an intern at LinkedIn during summer 2022. This work is not a part of his MIT research. Rahul Mazumder contributed to this work while he was a consultant for LinkedIn (in compliance with MIT’s outside professional activities policies).  This work is not a part of his MIT research.

\newpage
\bibliography{opt2022}
    \newpage
\appendix
\section{Proof of Lemma~\ref{lemma: remove-one-constraint}}

(i) Actually we can prove that both \eqref{LP-with-diversity} and \eqref{LP-with-diversity-ub} have the same set of optimal solutions as the following problem:
\begin{equation}\label{LP-with-diversity-ub-eq}
		\max_{X\in \R^{m\times n}} ~ c^\T X w
		\quad ~
		{\rm s.t.} ~  X \in \cS_{m,n}, ~ 
		 a^\T X w = b_2  . 
	\end{equation}
	Note that for any $X\in \cS_{m,n}$ with $a^\T X w <b_2$, one can write $\td X\in \cS_{m,n} $ as a convex combination of $X$ and some $X^0 \in \cX_{0}$, such that $a^\T \td X w = b_2 $ and $c^\T \td X w > c^\T X w $. 
	
(ii) can be proved via a similar argument as the proof of (i); 

(iii) holds true trivially.

\section{Proof of Lemma~\ref{lemma: screening}}

	Note that there exists $\al \in [0,1]$ such that $\lam^* = (1-\al) \lam_{\min} + \al \lam_{\max}$. Hence, 
	\begin{equation}
		\begin{aligned}
			(c-\lam^* a)_{i} &= (1-\al) ( c-\lam_{\min} a )_i + 
			\al ( c-\lam_{\max} a )_i \\
			&< 
			(1-\al) ( c-\lam_{\min} a )_{i_k} + 
			\al ( c-\lam_{\max} a )_{i_k} 
			= (c - \lam^* a)_{i_k}
			\nonumber
		\end{aligned}
	\end{equation}
for all $k \in [n]$, where the inequality is by \eqref{ineq-min-max}. 
As a result, we know that $i$ is not the top $n$ coordinates of $c - \lam^* a$, and hence the $i$-th row of the solution of \eqref{LP-with-diversity-ub} is zero.

\section{Algorithm Details and Discussions}

\begin{algorithm}[H]
	\caption{Bisection + tracing to solve \eqref{dual-problem}}
	\label{alg: bisection}
	\begin{algorithmic}
		\STATE \textbf{Input}: Tolerances $ \Delta > \dt \ge 0$.
		
		\STATE
		Initially set $\lam_{\min} = 0$, $\lam_{\max} = \infty$ and $ \lam = 1$. 
		
		\STATE \textbf{while} $\lam_{\max} - \lam_{\min}>\dt$:

		\STATE \quad Solve \eqref{def:g} and get $\cX_{\lam}$. 
		Compute $g_+'(\lam) $ and $ g_-'(\lam) $ using \eqref{grad-g}. 
		\STATE \quad \textbf{If} $g_+'(\lam) >0$: 
		
		\STATE \qquad If $\lam_{\max}- \lam_{\min} < \Dt $: 
		\STATE \qquad \quad Compute $K_R(\lam)$ and $g_+'(K_R(\lam))$ and $g_-'(K_R(\lam))$. 
		\STATE \qquad \quad If $g_+'(K_R(\lam)) \le 0 \le g_-'(K_R(\lam))$, \textbf{break} and \textbf{output} $\lam^* = K_R(\lam)$.

		\STATE \qquad If $\lam_{\max}<\infty$, then set $\lam_{\min} = \lam$, and $\lam = (\lam_{\max}+\lam_{\min})/2$; 
		\STATE \qquad Else, $\lam_{\max}=\infty$, then set $\lam_{\min} = \lam$ and $\lam = 2\lam$. 
		
		\STATE \quad \textbf{If} $g_-'(\lam) <0$: 
		
		\STATE \qquad If $\lam_{\max}- \lam_{\min} < \Dt $: 
		\STATE \qquad \quad Compute $K_L(\lam)$ and $g_+'(K_L(\lam))$ and $g_-'(K_L(\lam))$. 
		\STATE \qquad \quad If $g_+'(K_L(\lam)) \le 0 \le g_-'(K_L(\lam))$, \textbf{break} and \textbf{output} $\lam^* = K_L(\lam)$. 
		
		\STATE \qquad Set $\lam_{\max} = \lam$ and $\lam = (\lam_{\max}+\lam_{\min})/2$. 
		
		\STATE \quad \textbf{If} $ g_+'(\lam) \le 0 \le g_-'(\lam) $: 
		\STATE \qquad \textbf{output} $\lam^* = \lam$ as an optimal solution of \eqref{dual-problem}. 
		
		\STATE \textbf{end while}
		
		\STATE \textbf{Output} The interval $[\lam_{\min}, \lam_{\max}]$. 
		

	\end{algorithmic}
\end{algorithm}

\begin{algorithm}[H]
	\caption{Bisection + tracing + screening to solve \eqref{dual-problem}}
	\label{alg: bisection+screening}
	\begin{algorithmic}
		\STATE \textbf{Input}: Tolerances $ \Delta > \dt >0$.
		
		\STATE
		Initially set $\lam_{\min} = 0$, $\lam_{\max} = \infty$ and $ \lam = 1$. 
		 Set $\cI = [m]$. 
		
		\STATE \textbf{while} $\lam_{\max} - \lam_{\min}>\dt$:

		\STATE \quad Solve \eqref{def:g} and get $\cX_{\lam}$ (with $\cI$ in place of $[m]$). 
		
		\quad Compute $g_+'(\lam) $ and $ g_-'(\lam) $ using \eqref{grad-g}. 
		
		\STATE \quad \textbf{If} $g_+'(\lam) >0$: 
		
		\STATE \qquad If $\lam_{\max}- \lam_{\min} < \Dt $: 
		\STATE \qquad \quad Compute $K_R(\lam)$ and $g_+'(K_R(\lam))$ and $g_-'(K_R(\lam))$ (with $\cI$ in place of $[m]$). 
		\STATE \qquad \quad If $g_+'(K_R(\lam)) \le 0 \le g_-'(K_R(\lam))$, \textbf{break} and \textbf{output} $\lam^* = K_R(\lam)$.

		\STATE \qquad If $\lam_{\max}<\infty$, then set $\lam_{\min} = \lam$, and $\lam = (\lam_{\max}+\lam_{\min})/2$; 
		\STATE \qquad Else, $\lam_{\max}=\infty$, then set $\lam_{\min} = \lam$ and $\lam = 2\lam$. 
		
		\STATE \quad \textbf{If} $g_-'(\lam) <0$: 
		
		\STATE \qquad If $\lam_{\max}- \lam_{\min} < \Dt $: 
		\STATE \qquad \quad Compute $K_L(\lam)$ and $g_+'(K_L(\lam))$ and $g_-'(K_L(\lam))$ (with $\cI$ in place of $[m]$). 
		\STATE \qquad \quad If $g_+'(K_L(\lam)) \le 0 \le g_-'(K_L(\lam))$, \textbf{break} and \textbf{output} $\lam^* = K_L(\lam)$. 
		
		\STATE \qquad Set $\lam_{\max} = \lam$ and $\lam = (\lam_{\max}+\lam_{\min})/2$. 
		
		\STATE \quad \textbf{If} $ g_+'(\lam) \le 0 \le g_-'(\lam) $: 
		\STATE \qquad \textbf{output} $\lam^* = \lam$ as an optimal solution of \eqref{dual-problem}. 
		
		\STATE \quad (Screening to shrink $\cI$)
		\STATE \quad \textbf{If} $\lam_{\max}<\infty$ :
			
			\STATE \qquad
			1. Let $i_1,...,i_n$ be the indices in $\cI$ corresponding to top $n$ elements of $(c- \lam a)$.
			
			\qquad 
			(Break ties arbitrarily, if any). 
			
			\qquad 2. Find out 
			\begin{equation}
				\small
				\cI_{drop} := \Big\{
				i\in \cI ~\Big|~ \text{Inequalities in \eqref{ineq-min-max} are satisfied}
				\Big\} \nonumber
			\end{equation}
			
			\qquad 
			and update $\cI = \cI \setminus \cI_{drop}$.

		\STATE \textbf{end while}
		
		\STATE \textbf{Output} The interval $[\lam_{\min}, \lam_{\max}]$. 
		
	\end{algorithmic}
\end{algorithm}

Note that in Algorithm~\ref{alg: bisection+screening}, a set of active indices $\cI$ is maintained and updated in each iteration. For the computations of $\cX_{\lam}$ and $K_{L}(\lam)$ (or $K_R(\lam)$), only indices in $\cI$ are considered, which will have computational costs $O(|\cI| \log(|\cI|))$ and $|\cI| n$ respectively. If the screening step can discard a large set of candidates in the first few iterations of the Algorithm~\ref{alg: bisection+screening}, then the overall cost of Algorithm~\ref{alg: bisection+screening} is significantly lower than that of Algorithm~\ref{alg: bisection}.




\section{Recovering primal solutions}\label{section: recovering primal solutions}

If Algorithm~\ref{alg: bisection} finds the exact dual optimal solution $\lam^*$, the (primal) optimal solution of \eqref{LP-with-diversity-ub} can be recovered directly. Note that Algorithm~\ref{alg: bisection} also makes $\cX_{\lam^*}$ available. Since $\lam^*$ is the dual optimal solution, $g_-'(\lam^*) \le 0$ and $g_+'(\lam^*) \ge 0$, implying $\min_{X\in \cX_{\lam^*}} \{ a^\T X w \} \le 
	 b_2 \le  \max_{X\in \cX_{\lam^*}} \{ a^\T X w \}\ayan{.} $
Let $X^1$ and $X^2$ be two extreme points of $\cX_{\lam^*}$ such that $a^\T X^1 w \le b_2 $ and $a^\T X^2 w \ge b_2 $. Define 
$\rho := (b_2 - a^\T X^2 w) /( a^\T X^1 w - a^\T X^2 w  )$ and 
$X^* := \rho X^1 + (1-\rho) X^2, $
then $X^* \in \cS_{m,n}$ and $a^\T X^* w = b_2$. By the KKT condition, $X^*$ is an optimal solution of problem \eqref{LP-with-diversity-ub}. 

\section{Data generation}\label{section: data generation}

We set the weights $w$ as $w_j = 1/\log(1+j)$ for all $j\in [n]$ as commonly used in the Discounted Cumulative Gain (DCG) measure \cite{singh2018fairness}. We generate vector pairs $(a, c)\in \R^{m\times 2}$ such that $\{ (a_i, c_i ) \}_{i=1}^n$ are i.i.d. bi-variate Gaussian random variable with mean $\Ex (a_i) = \Ex(c_i) = 0$, variance
$\var(a_i) = \var(c_i) = 1 $ and covariance ${\rm cov}(a_i, c_i) = \alpha$ for some chosen $\alpha \in (0,1)$. Note that $\alpha$ controls the angle between the vectors $a$ and $c$; for a larger value of $ \alpha $, the angle between $a$ and $c$ are smaller (with high probability). We set $\al = 0.5$. 
Finally, to generate $b_1$ and $b_2$, we compute $\bar X \in \argmax_{X\in \cS_{m,n}} \{c^\T X w\}$, and evaluate $a^\T \bar X w$; if $a^\T \bar X w >0 $ (which happens with high probability if $\al$ is not too small), take 
$b_2 = -b_1 = 0.8 (a^\T \bar X w)$; otherwise drop this instance and regenerate one with $a^\T \bar X w >0 $. With these values of $b_1$ and $b_2$, the constraint $b_1 \le a^\T X w \le b_2$ must be active at the optimal solution of \eqref{LP-with-diversity}.

\section{Performance on Real Datasets}
We evaluate the performance of Algorithm~\ref{alg: bisection+screening} on proprietary real-world datasets that arise in a recommender system of people recommendations to generate connections. To improve the diversity across different groups of the recommended people, we encode the group a candidate $i$ belongs to as $a_i \in \{-1,0,1\}$. Here $a_i = 1$ represents group $1$, $a_i = -1$ stands for group $2$, and $a_i = 0$ denotes other groups. Information about the specific feature according to which such groups are formed are withheld here due to privacy considerations. We choose some $b>0$ and set $b_2 = -b_1 = b$ in problem \eqref{LP-with-diversity}. 

We focus on two instances of the problem with $n = 10$ slots and $m = 500$ or $1000$ candidates respectively. To evaluate the robustness of Algorithm~\ref{alg: bisection+screening}, we created $20$ noised replications of problem \eqref{LP-with-diversity} with noised scores $c_i$ ($i\in [m]$). More precisely, we set $c = c^0 + 0.2\frac{\|c^0\|_2}{\|\ep\|_2} \ep$, where $c^0\in \R^m$ is the original utility scores, and $\ep \in \R^m$ is a random vector with i.i.d. $N(0,1)$ coordinates. Other data ($a$ and $w$) assume the same values as the original data. 


\begin{figure}[H] 
	\centering 
	\includegraphics[height=5cm,width=8cm]{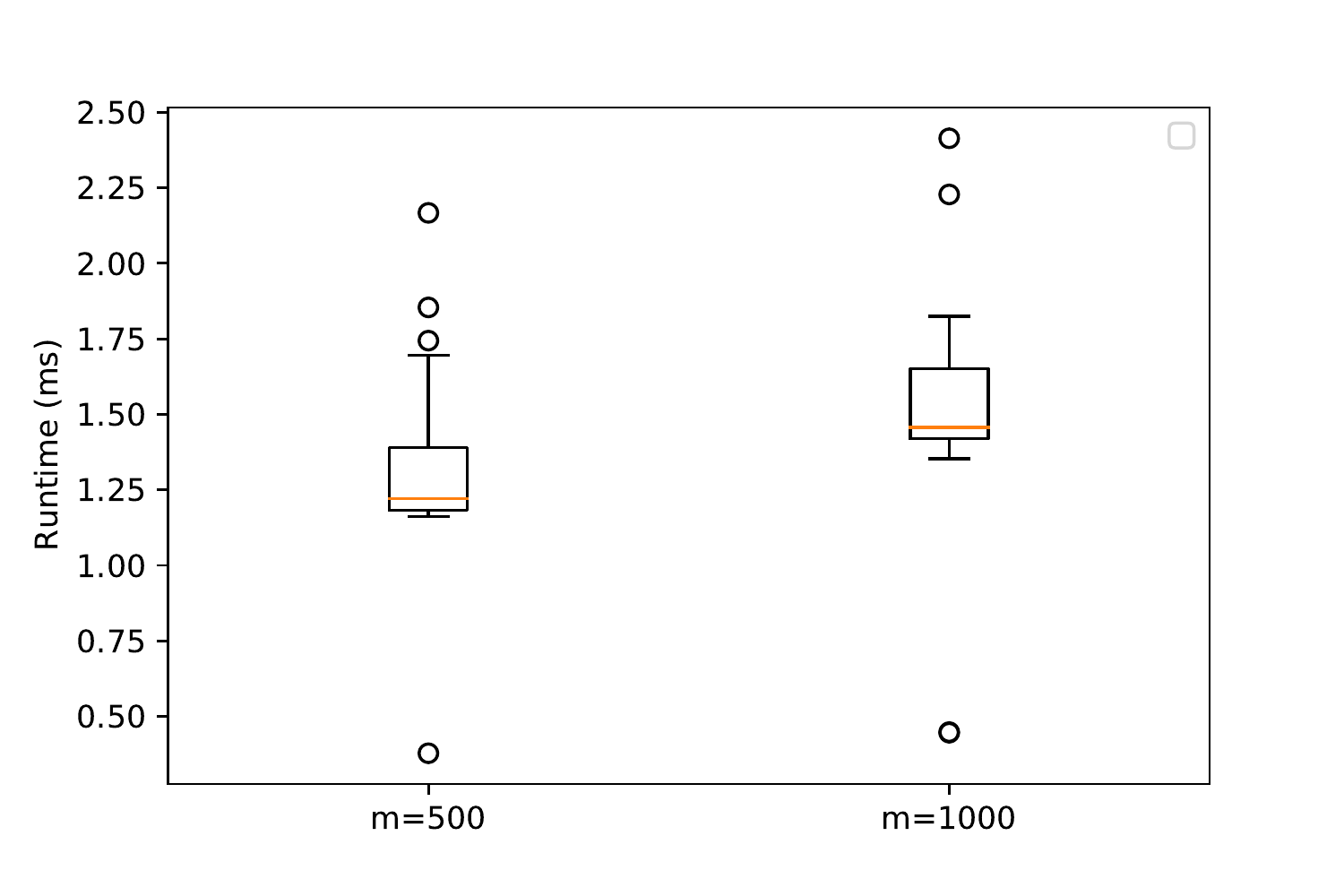}
	\caption{Performance of Algorithm~\ref{alg: bisection+screening} on real datasets
	} 
	\label{figure: real}
\end{figure}

Figure~\ref{figure: real} presents the box plot of Algorithm~\ref{alg: bisection+screening} runtimes for these 20 replications. For $m = 500$, most replications have a runtime between $1$ ms and $2.25$ ms; there are one or two replications that have a runtime below $0.5$ ms -- these replications correspond to the case when the diversity constraint $b_1 \le a^\T X w \le b_2$ is not active at the optimal solution, and Algorithm~\ref{alg: bisection+screening} terminates earlier. A similar argument applies to the examples with $m = 1000$. As a conclusion, Algorithm~\ref{alg: bisection+screening} is quite robust across different runs with noised data.

\section{Real-world Use Cases}\label{sect:real-world-use-case}

Below, we present a few concrete real-world use cases which can potentially leverage and benefit from our specialized solver.

\begin{itemize}
\item One might be interested in providing more diverse exposure to different types of content creators on a social network platform. For example, the platform might want to show diverse video recommendations from creators from different countries or states.

\item Introducing systematic diversity in movie recommendations can 
lead to better exploration and user retention. For instance, the recommendations can span across genres such as ``romance", ``action", ``sci-fi", ``comedy", etc., to better exploit the interests of the users. 

\item Our algorithms can also find applications related to ride-sharing apps. Experienced and frequent drivers on such platforms usually have enough exposure to riders. However, more balanced exposure to new or less frequent drivers may help ride-sharing platforms achieve more sustainable growth. Of course, the same strategy may also find relevance in other two-sided marketplace platforms, such as house sharing, restaurant recommendations, etc.
\end{itemize}

%
%



\end{document}